\documentclass[10pt]{amsart}

\usepackage{myStyle}
\graphicspath{{./Figs/}}

\newtheorem*{theorem}{Theorem}

\begin{document}

\title{A Short Note on Mapping Cylinders}

\author{Alex Aguado}
\address{Department of Mathematics, Duke University, Durham, NC 27708-0320}
\email{aaguado@math.duke.edu}



\begin{abstract}
Given a homotopy equivalence $f:X \absRightarrow Y$ we assemble well known pieces and unfold them into an explicit formula for a strong deformation retraction of the mapping cylinder of $f$ onto its top.
\end{abstract}

\maketitle

\section{Setup}

The mapping cylinder $M_{f}$ of a map $f:X\rightarrow Y$ between two topological spaces is the quotient space of the disjoint union $X\times I+Y$ by the equivalence relation arising from the identifications $(x,0)\sim f(x)$. The quotient map $q:X\times I+Y\rightarrow M_{f}$ embeds $X\times\{1\}$ and $Y$ as closed subsets of $M_{f}$. Let $\widetilde{X}=q(X\times\{1\})$, and $\widetilde{Y}=q(Y)$. We will denote by $[x,t]$ and $[y]$ the equivalence classes of $(x,t)\in X\times I$, and $y\in Y$, respectively.

\begin{theorem} A map $f:X\rightarrow Y$ between two topological spaces is a homotopy equivalence iff $\widetilde{X}$ is a strong deformation retract of $M_{f}$.
\end{theorem}

This theorem is well known, and its origins can be traced back to the results presented by Fox \cite{Fox} and Fuchs \cite{Fuchs}. Modern proofs are usually presented within the contexts of cofibrations or the homotopy extension property (e.g. Corollary 0.21 in \cite{Hatcher}).

The backward implication is usually dealt with by giving an explicit formula for a strong deformation retraction of $M_{f}$ onto its bottom $\widetilde{Y}$ (e.g. the map $H_{1}$ given below).

The forward implication has a conceptual proof that we believe is better than a proof by a formula. However, we also believe that it is desirable to have a formula because of computer modeling, which has become an important research tool across the scientific community. 

In particular, the computational aspects of the problem treated in this note may be resolved by breaking it down into its fundamental pieces, implementing each separately on a computer, and composing them by function calls. On the other hand, one may assemble the same fundamental pieces and unfold them into a formula, which is easily implemented on a computer. To the best of out knowledge the first method was the only choice. We provide the second. Deciding which one is computationally more efficient is a task for expert programmers.	

\section{Tools}

Suppose that $f:X\rightarrow Y$ is a homotopy equivalence with homotopy inverse $g:Y\rightarrow X$. Then we have homotopies $F:X\times I\rightarrow X$ from $g\circ f$ to $\bm{1}_{X}$, and $G:Y\times I\rightarrow Y$ from $f\circ g$ to $\bm{1}_{Y}$.

We want to construct a strong deformation retraction of $M_{f}$ onto its top $\widetilde{X}$ in terms of $f$, $g$, $F$ and $G$. The following are the tools we will need:

\subsection{A deformation retraction of $\bm{M_{f}}$ onto its top $\bm{\widetilde{X}}$:}\label{sectDefRet}\medskip
\begin{enumerate}
\item Define a map $H_{1}:M_{f}\times I\rightarrow M_{f}$ by\medskip
\[
\bigl([x,t],s\bigr)\longmapsto [x,t(1-s)]\;\;\text{and}\;\;\bigl([y],s\bigr)\longmapsto [y]\medskip
\]
This is a well defined strong deformation retraction from $\bm{1}_{M_{f}}$ to the retraction $r:M_{f}\rightarrow M_{f}$ given by the canonical projection of the mapping cylinder onto its bottom $\widetilde{Y}$. That is, $r$ is defined by $[x,t]\mapsto [x,0]$, and $[y]\mapsto [y]$.\medskip
\item Use $f$ and $G$ to define a map $H_{2}:M_{f}\times I\rightarrow M_{f}$ by\medskip
\[
\bigl([x,t],s\bigr)\longmapsto [G(f(x),1-s)]\;\; \text{and}\;\;\bigl([y],s\bigr)\longmapsto [G(y,1-s)]\medskip
\]
This is a well defined homotopy from $r$ to a map $h:M_{f}\rightarrow M_{f}$ defined by $[x,t]\mapsto[g\circ f(x),0]$ and $[y]\mapsto [g(y),0]$.\medskip
\item Use $g$ and $F$ to define a map $H_{3}:M_{f}\times I\rightarrow M_{f}$ by\medskip
\[
\bigl([x,t],s\bigr)\longmapsto[F(x,st),s]\;\;\text{and}\;\;\bigl([y],s\bigr)
\longmapsto [g(y),s]\medskip
\]
This is a well defined homotopy from the map $h$ to a retraction $r':M_{f}\rightarrow M_{f}$ of the mapping cylinder onto its top $\widetilde{X}$ defined by\medskip
 \[
 [x,t]\longmapsto [F(x,t),1]\;\;\text{and}\;\;[y]\longmapsto [g(y),1].\medskip
 \]
\end{enumerate}

The concatenation $H=H_{1}*H_{2}*H_{3}$ is the desired deformation retraction from $\bm{1}_{M_{f}}$ to a retraction $r'$ of $M_{f}$ onto its top.

\subsection{The Homotopy extension property of $\bm{(M_{f}\times I,\widetilde{X}\times I)}$ made explicit:} We need to construct a retraction $R:M_{f}\times I\times I\rightarrow M_{f}\times I\times I$ of $M_{f}\times I\times I$ onto $M_{f}\times I\times\{0\}\cup \widetilde{X}\times I\times I$. This is achieved by defining first a retraction $\varphi:I^{2}\rightarrow I^{2}$ of $I^{2}$ onto $I\times\{0\}\cup\{1\}\times I$ via radial projection from the point $(0,2)$, as the figure below illustrates.\hfill
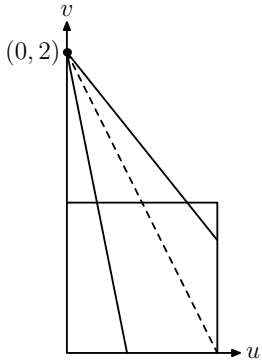
\smallskip
\begin{wrapfigure}[6]{l}[-1cm]{3.5cm}
\includegraphics[scale=0.8]{figures.1}
\end{wrapfigure}
\[
\hspace{4.5cm}\varphi(u,v)=\begin{cases}
              \left(\dfrac{2u}{2-v},0\right)& \text{if\, $v\leq 2-2u$}\\\\
              \left(1,\dfrac{2u+v-2}{u}\right)& \text{if\, $v\geq 2-2u$.}
\end{cases}\vspace{0.8cm}
\]
If we let $\varphi_{1}$ and $\varphi_{2}$ denote the components of $\varphi$ the retraction $R$ is defined by
\[
\hspace{4.5cm}\bigl([x,t],s,l\bigr)\longmapsto \bigl([x,\varphi_{1}(t,l)],s,\varphi_{2}(t,l)\bigr)\;\;\;\mbox{and}\;\;\;
\bigl([y],s,l\bigr)\longmapsto\bigl([y],s,0\bigr).\vspace{0.5cm}
\]

\section{Construction}

In order to facilitate clarity a point in $M_{f}$ will be denoted by the letter $p$. If the point is in $\widetilde{X}$ (i.e. $p=[x,1]$) we will denote it by $\tilde{p}$.

Our goal is to modify the deformation retraction $H$ constructed earlier so that it leaves $\widetilde{X}$ fixed.

We begin by defining a homotopy $K:M_{f}\times I\rightarrow M_{f}$ by\medskip
\[
K(p,s)=\begin{cases}
        H^{-1}(p,1-2s)& \text{if\, $0\leq s\leq \dfrac{1}{2}$}\\\\
        H^{-1}(r'(p),2s-1)&     \text{if\, $\dfrac{1}{2}\leq s\leq 1$,}
\end{cases}\medskip
\]
where $H^{-1}$ denotes the \emph{inverse} of $H$ (i.e. $H$ with $s$ running backward). Roughly speaking, the first part of $K$ is simply $H^{-1}$ twice as fast, while the second half is a map that begins and ends with $r'$, and is homotopic\footnote{Because $\bm{1}_{M_{f}}\sim r'$ implies that $\bm{1}_{M_{f}}\times\bm{1}_{I}\sim r'\times\bm{1}_{I}$, and therefore $H^{-1}\sim H^{-1}\circ (r'\times\bm{1}_{I})$.} to $H^{-1}$. It easy to check that $K$ is a well defined homotopy from $\bm{1}_{M_{f}}$ to $r'$, and that its restriction to $\widetilde{X}\times I$ is a homotopy beginning and ending with $\bm{1}_{\widetilde{X}}$. Moreover, it behaves well with respect to time reversal in the sense that $K(\tilde{p},1-s)=K(\tilde{p},s)$.

The next step will be to define a homotopy $L:(\widetilde{X}\times I)\times I\rightarrow M_{f}$. We use $K$ and the diagram below to accomplish this.\hfill\phantom{}
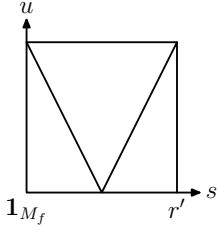
\begin{wrapfigure}[4]{l}[-1cm]{3cm}
\includegraphics[scale=0.8]{figures.3}
\end{wrapfigure}
Roughly speaking, we want $L(\tilde{p},s,u)$ to be independent of $u$ below the `V', and  independent of $s$ above the `V'. More formally:\bigskip
\[
\hspace{3.5cm}L(\tilde{p},s,u)=\begin{cases}
K(\tilde{p},s)& \text{if\, $u\leq|2s-1|$}\\\\
K\left(\tilde{p},\dfrac{1-u}{2}\right)& \text{if\, $2s-u\leq 1\leq 2s+u$.}
\end{cases}\vspace{0.5cm}
\]
This is a well defined\footnote{Because $K(\tilde{p},1-s)=K(\tilde{p},s)$ implies that $K(\tilde{p},\frac{1-u}{2})=K(\tilde{p},\frac{1+u}{2})$.} homotopy from $K\left|_{\widetilde{X}\times I}\right.$ to a map sending $(\tilde{p},s)$ to $\tilde{p}$. Moreover, $L(\tilde{p},s,u)=\tilde{p}$ for $(s,u)\in \partial I\times I\cup I\times \{1\}$. Therefore, looking at the diagram above we would like to extend $L$ to the whole mapping cylinder and then follow this extension from the lower left corner $\bm{1}_{M_{f}}$ along the left, top and right edges to the right corner $r'$. Indeed, define a map $K_{0}:M_{f}\times I\times \{0\}\rightarrow M_{f}$ by $K_{0}(p,s,0)=K(p,s)$, and combine it with $L$ to give a map\footnote{This map is unambiguously defined because $K_{0}=L$ on the overlap $\widetilde{X}\times I\times \{0\}$. Continuity follows from the gluing lemma because $\widetilde{X}$ is a closed subset of in $M_{f}$.} $\left(K_{0},L\right):M_{f}\times I\times \{0\}\cup \widetilde{X}\times I\times I\rightarrow M_{f}$. Then, define $L'$ by the composition\medskip
\[
\begin{psmatrix}
M_{f}\times I\times I   &M_{f}\times I\times \{0\}\cup \widetilde{X}\times I\times I   &M_{f},
\end{psmatrix}\medskip
\everypsbox{\scriptstyle}
\psset{shortput=nab,linewidth=0.56pt,nodesep=3pt,arrowsize=4pt,arrows=->}
\ncline{1,1}{1,2}^[npos=0.5]{R}
\ncline{1,2}{1,3}^[npos=0.5]{\left(K_{0},L\right)}
\]
This map clearly extends $L$ and satisfies $L'(p,s,0)=K(p,s)$. Finally, define a map $\Gamma:M_{f}\times I \rightarrow M_{f}$ by\bigskip
\[
\Gamma(p,s)=\begin{cases}
L'(p,0,3s)&     \text{if\, $0\leq s\leq\dfrac{1}{3}$}\\\\
L'(p,3s-1,1)&   \text{if\, $\dfrac{1}{3}\leq s\leq\dfrac{2}{3}$}\\\\
L'(p,1,3-3s)&   \text{if\, $\dfrac{2}{3}\leq s\leq 1$}.
\end{cases}\bigskip
\]
It is easy to check that this map verifies the following:\medskip
\begin{itemize}
\item $\Gamma(p,0)=L'(p,0,0)=K(p,0)=p$,\smallskip
\item $\Gamma(p,1)=L'(p,1,0)=K(p,1)=r'(p)$,\smallskip
\item $\Gamma(\tilde{p},s)=\tilde{p}$ because $L'$ is an extension of $L$.
\end{itemize}\medskip
Therefore, $\Gamma$ is the desired strong deformation retraction from the identity $\bm{1}_{M_{f}}$ to a retraction $r'$ of the mapping cylinder onto its top $\widetilde{X}$.

\section{Formula}

If we unravel $\Gamma$ we will obtain the desired formula in terms of $f$, $g$, $F$ and $G$. It is given in the next two pages.\smallskip
\newpage
\[\hspace{-0.80cm}
\Gamma\bigl([x,t],s\bigr)=
\begin{cases}
\text{if\, $0\leq s\leq \dfrac{1}{3}$} &\begin{cases}
\left[x,\dfrac{2t}{2-3s}\right]& \text{if\, $s\le\dfrac{2-2t}{3}$}\\\\
[x,1]& \text{if\, $s\geq\dfrac{2-2t}{3}$}
\end{cases}\\[2cm]
\text{if\, $\dfrac{1}{3}\leq s\leq \dfrac{2}{3}$}& \begin{cases}
\text{if\, $t\leq\dfrac{1}{2}$}&
\begin{cases}
\left[x,2t(7-18s)\right]& \text{if\, $\dfrac{1}{3}\leq s\leq \dfrac{7}{18}$}\\[0.4cm]
\left[G(f(x),8-18s)\right]& \text{if\, $\dfrac{7}{18}\leq s\leq \dfrac{8}{18}$}\\[0.4cm]
\left[F\left(x,2t(18s-8)\right),18s-8\right]& \text{if\, $\dfrac{8}{18}\leq s\leq \dfrac{1}{2}$}\\[0.4cm]
\left[F\left(F(x,2t),10-18s\right),10-18s\right]& \text{if\, $\dfrac{1}{2}\leq s\leq \dfrac{10}{18}$}\\[0.4cm]
\left[G\left(f\circ F(x,2t),18s-10\right)\right]& \text{if\, $\dfrac{10}{18}\leq s\leq \dfrac{11}{18}$}\\[0.4cm]
\left[F(x,2t),18s-11\right]& \text{if\, $\dfrac{11}{18}\leq s\leq \dfrac{2}{3}$}
\end{cases}\\[3.5cm]
\text{if\, $t\geq\dfrac{1}{2}$}&
\begin{cases}
\left[x,7-18s\right]\hphantom{OOOOOOOOOOO||}& \text{if\, $\dfrac{1}{3}\leq s\leq \dfrac{7}{18}$}\\[0.4cm]
\left[G(f(x),8-18s)\right]& \text{if\, $\dfrac{7}{18}\leq s\leq \dfrac{8}{18}$}\\[0.4cm]
\left[F\left(x,18s-8\right),18s-8\right]& \text{if\, $\dfrac{8}{18}\leq s\leq \dfrac{1}{2}$}\\[0.4cm]
\left[F\left(x,10-18s\right),10-18s\right]& \text{if\, $\dfrac{1}{2}\leq s\leq \dfrac{10}{18}$}\\[0.4cm]
\left[G(f(x),18s-10)\right]& \text{if\, $\dfrac{10}{18}\leq s\leq \dfrac{11}{18}$}\\[0.4cm]
\left[x,18s-11\right]& \text{if\, $\dfrac{11}{18}\leq s\leq \dfrac{2}{3}$}
\end{cases}
\end{cases}\\[7.65cm]
\text{if\, $\dfrac{2}{3}\leq s\leq 1$}& \begin{cases}
\left[F\left(x,\dfrac{2t}{3s-1}\right),1\right]& \text{if\, $s\geq \dfrac{1+2t}{3}$}\\\\
\left[x,1\right]& \text{if\, $s\leq \dfrac{1+2t}{3}$}
\end{cases}
\end{cases}
\]
\newpage
\[
\Gamma\bigl([y],s\bigr)=\begin{cases}
[y]& \text{if\, $0\leq s\le\dfrac{7}{18}$}\\[0.4cm]
\left[G(y,8-18s)\right]& \text{if\, $\dfrac{7}{18}\leq s\le\dfrac{8}{18}$}\\[0.4cm]
[g(y),18s-8]& \text{if\, $\dfrac{8}{18}\leq s\le\dfrac{1}{2}$}\\[0.4cm]
\left[F\left(g(y),10-18s\right)\right]& \text{if\, $\dfrac{1}{2}\leq s\le\dfrac{10}{18}$}\\[0.4cm]
\left[G\left(f\circ g(y),18s-10\right)\right]& \text{if\, $\dfrac{10}{18}\leq s\le\dfrac{11}{18}$}\\[0.4cm]
[g(y),18s-11]& \text{if\, $\dfrac{11}{18}\leq s\le\dfrac{2}{3}$}\\[0.4cm]
[g(y),1]& \text{if\, $\dfrac{2}{3}\leq s\le 1$}
\end{cases}
\]

\end{document}